\documentclass[11pt]{amsart}

\usepackage{amsmath,amssymb,amscd,amsthm,graphicx,enumerate}

\usepackage{hyperref}
\hypersetup{breaklinks=true}
\usepackage{graphicx,float}

\theoremstyle{plain}

\newtheorem{theorem}{Theorem}[section]

\newtheorem{proposition}[theorem]{Proposition}

\newtheorem{corollary}[theorem]{Corollary}
\newtheorem{conjecture}[theorem]{Conjecture}
\newtheorem{problem}[theorem]{Problem}

\theoremstyle{definition}
\newtheorem{definition}[theorem]{Definition}

\newtheorem{example}[theorem]{Example}

\begin{document}

\title[Mean curvature flow through singularities]{Mean curvature flow through singularities}
\author{Robert Haslhofer}

\date{\today}

\maketitle

\begin{abstract} We first give a general introduction to the mean curvature flow, and then discuss fundamental results established over the last 10 years that yield a precise theory for the flow through singularities in $\mathbb{R}^3$. With the aim of developing a satisfying theory in higher dimensions, we then describe our recent classification of all noncollapsed singularities in $\mathbb{R}^4$. Finally, we provide a detailed discussion of open problems and conjectures.
\end{abstract}

\section{Introduction}
Mean curvature flow describes surfaces that evolve to decrease their area as efficiently as possible. More precisely, a smooth family of embedded surfaces $\{M_t\subset \mathbb{R}^{3}\}_{t\in I}$ moves by mean curvature flow if
\begin{equation}\label{eq_mcf}
\partial_t x = \vec{H}(x)
\end{equation}
for all $x\in M_t$ and $t\in I$. Here, $I\subset \mathbb{R}$ is an interval, $\partial_t x$ is the normal velocity at $x$, and $\vec{H}(x)$ is the mean curvature vector at $x$. Concretely, if we write $\vec{H}=H\nu$, where $\nu$ is the inwards unit normal, then $H$ is given by the sum of the principal curvatures, namely $H=\kappa_1+\kappa_2$. The simplest example of a mean curvature flow is a family of round shrinking spheres $M_t=S^2({r(t)})$, where $r(t)=\sqrt{R^2-4t}$. However, most mean curvature flows are not given by explicit formulas, and one rather seeks for a quantitative description of the evolution.

The mean curvature flow first appeared as a model for evolving interfaces in material science \cite{Mullins}, and its mathematical study was pioneered by Brakke \cite{brakke} and Huisken \cite{Huisken_convex}. It is arguably the most natural evolution equation in extrinsic geometry, a close cousin to Hamilton's Ricci flow in intrinsic geometry \cite{Hamilton_survey}, and has been a very active field of research in pure and applied mathematics over the last 50 years.

To understand the flow from a PDE perspective, recall that given any point $p$ on a surface $M$, if we express the surface locally as a graph of a function $u$ over the tangent space $T_pM$, then the principal curvatures $\kappa_1(p)$ and $\kappa_2(p)$ are the eigenvalues of $\mathrm{Hess}(u)(p)$. Hence, \eqref{eq_mcf} is a (nonlinear) heat-type equation, and one hopes that the flow will deform the initial surface into a nicer one. This indeed works very well:

\begin{theorem}[Huisken's convergence theorem \cite{Huisken_convex}]
Let $M_0\subset \mathbb{R}^{3}$ be a closed embedded surface. If $M_0$ is convex, then the mean curvature flow $\{M_t\}_{t\in[0,T)}$ starting at $M_0$ converges to a round point.
\end{theorem}

The convex case is of course very special.
In more general situations, one encounters local singularities:

\begin{example}[neckpinch singularity \cite{AltschulerAngenentGiga,AngenentVelazquez}]
If $M_0$ has the topology of a sphere but the geometry of a dumbbell, then the neck pinches off. As blowup limit (i.e. as limit of rescalings of the flow around the singularity along a sequence of rescaling factors $\lambda_i\to \infty$) we get a selfsimilarly shrinking round cylinder. There is also a degenerate variant of this example, where as blowup limit along suitable tip points one gets a selfsimilarly translating bowl.
\end{example}

In the study of mean curvature flow the main task is to analyze the formation of singularities:
\begin{itemize}
\item How do singularities look like?
\item Can we continue the flow through singularities?
\item What is the size and the structure of the singular set?
\item Is the evolution through singularities unique or nonunique?
\end{itemize}

This article\footnote{Some parts draw from my summer school lecture notes \cite{Haslhofer_summer_school}, which I recommend for a more leisurely introduction.} is organized as follows. In Section \ref{section2}, we review some basic notions and fundamental tools that are crucial for the analysis of singularities. In Section \ref{section3}, we discuss the by now highly developed theory of mean curvature flow through singularities in $\mathbb{R}^3$. In Section \ref{section4}, we describe our recent classification of all noncollapsed singularities in $\mathbb{R}^4$. To conclude, in Section \ref{section5}, we discuss several further directions and open problems.

\textbf{Acknowledgements.} I am very grateful to Kyeongsu Choi, Wenkui Du and Or Hershkovits for fruitful long-term collaborations, whose results are described here, and to Tom Ilmanen and Bruce Kleiner for getting me excited about mean curvature flow as a student and postdoc. My research has been supported by an NSERC Discovery Grant.

\section{Basic notions and fundamental tools}\label{section2}

In this section, we first review some basic properties and then discuss preserved quantities, Huisken's monotonicity formula, $\varepsilon$-regularity, and weak solutions.

\subsection{Basic properties}

Let us collect a few basic properties of mean curvature flow that will play an important role throughout:

\begin{itemize}
\item \textbf{Short time existence:} Given any smooth compact initial surface $M_0\subset \mathbb{R}^{3}$, there exists a unique smooth solution $\{M_t\}_{t\in [0,T)}$ of \eqref{eq_mcf} starting at $M_0$. The maximal existence time $T$ is characterized by blow up of the second fundamental form, namely $\lim_{t\to T}\max_{M_t}|A|=\infty$.

\item \textbf{Avoidance principle:} If $\{M_t\}$ and $\{N_t\}$ are two  mean curvature flows, say at least one of them compact, then $\mathrm{dist}(M_t,N_t)$ is nondecreasing in time. In particular, if $M_{t_0}$ and $N_{t_0}$ are disjoint, then so are $M_t$ and $N_t$ for all $t\geq t_0$. Similarly, the flow does not bump into itself, i.e. embeddedness is also preserved.

\item \textbf{Parabolic scaling:} If $\mathcal{M}=\{M_t\}$ is a mean curvature flow, and $\mathcal{D}_\lambda \mathcal{M}$ denotes the flow obtained from $\mathcal{M}$ by the change of variables $(x,t)\mapsto (\lambda x, \lambda^2 t)$, then $\mathcal{D}_\lambda \mathcal{M}$ again solves \eqref{eq_mcf}.
\end{itemize}

Indeed, short time existence follows from standard parabolic theory (see e.g. \cite{HP,MM}), and the avoidance principle is just a geometric incarnation of the maximum principle. In particular, by the avoidance principle any compact surface becomes extinct in finite time, since we can engulf it by a sphere of some large radius $R<\infty$, which shrinks to a point at time $T=R^2/4$. Finally, the parabolic scaling,  while a trivial consequence of the chain rule, is of fundamental importance for the singularity analysis throughout: To analyze how the flow $\mathcal{M}$ looks like under the microscope near a space-time point $X_0=(x_0,t_0)$, one considers the flow $\mathcal{D}_\lambda (\mathcal{M}-X_0)$, which is obtained from $\mathcal{M}$ by shifting $X_0$ to the origin and parabolically rescaling by $\lambda$. E.g. if there is a neckpinch singularity at $X_0$, then for $\lambda$ large the flow $\mathcal{D}_\lambda (\mathcal{M}-X_0)$ looks almost like a round shrinking cylinder.

Another basic observation is that the evolution equation \eqref{eq_mcf} implies evolution equations for all induced geometric quantities. E.g. the evolution of the induced metric $g_{ij}=\partial_i X\cdot \partial_j X$ in some local parametrization $X(\cdot, t)$ of $M_t$, is obtained by computing
\begin{equation}
\partial_t g_{ij}=2\partial_i (H\nu)\cdot \partial_j X=2H\partial_i \nu\cdot \partial_j X=-2HA_{ij},
\end{equation}
where the first step used the mean curvature flow equation, and the last step used the definition of the second fundamental form.
In particular, remembering that the derivative of the determinant is given by the trace, one sees that under mean curvature flow the area element $d\mu=\sqrt{\det g_{ij}} \, d^2x$ evolves by $\partial_td\mu=-H^2 d\mu$, hence
\begin{equation}\label{eq_areamon}
\frac{d}{dt}\textrm{Area}(M_t)=-\int_{M_t} H^2 d\mu.
\end{equation}
This formula captures the important fact that the mean curvature flow is the gradient flow of the area functional, i.e. indeed the most efficient evolution to decrease area. Finally, another standard computation yields
\begin{equation}\label{ev_mean_curv}
\partial_t H=\Delta H+|A|^2 H,
\end{equation}
where $|A|^2=\kappa_1^2+\kappa_2^2$ denotes the squared norm of the second fundamental form.

\subsection{Preserved quantities and convexity estimate} 
Like for any evolution equation, it is very useful to identify preserved quantities. First, in light of \eqref{ev_mean_curv}, the maximum principle immediately yields that mean-convexity, namely the inequality $H\geq 0$, is preserved under mean curvature flow:

\begin{proposition}[mean-convexity]\label{thm_andrews}
 If $M_0$ is mean-convex, then so is $M_t$ for all $t\geq 0$.
\end{proposition}

Mean-convexity is on the one hand flexible enough to allow for the formation of local singularities, e.g. neckpinches and degenerate neckpinches, and flexible enough for interesting applications,\footnote{The most famous application for inverse mean curvature flow is the proof of the Penrose inequality by Huisken-Ilmanen \cite{HuiskenIlmanen}. For some other applications of mean-convex flows see e.g. \cite{Schulze_isoperimetric,huisken-sinestrari3,BHH,HaslhoferKetover,WangZhou,HaslhoferKetover2,LiokumovichMaximo}.} but on the other hand rigid enough to obtain a detailed description of singularities.

The following quantitative notion of embeddedness plays a key role in the theory of mean-convex flows:

\begin{definition}[noncollapsing \cite{ShengWang,andrews1,HK1}]\label{def_andrews_static}
A closed embedded mean-convex surface $M$ is called \emph{$\alpha$-noncollapsed},
if each point $p\in M$ admits interior and exterior balls tangent at $p$ of radius $\alpha/H(p)$.
\end{definition}

A beautiful argument of Andrews shows:

\begin{theorem}[noncollapsing \cite{andrews1}]\label{thm_andrews}
 If $M_0$ is $\alpha$-noncollapsed, then so is $M_t$ for all $t\geq 0$.
\end{theorem}

Indeed, Andrews first observed that $\alpha$-noncollapsing is equivalent to the inequalities
\begin{equation}
 \sup_{y\neq x} \frac{2\langle X(y,t)-X(x,t),\nu(x,t)\rangle}{|X(y,t)-X(x,t)|^2}\leq \frac{H(x,t)}{\alpha}\quad\mathrm{and}\quad \inf_{y\neq x} \frac{2\langle X(y,t)-X(x,t),\nu(x,t)\rangle}{|X(y,t)-X(x,t)|^2}\geq -\frac{H(x,t)}{\alpha},
\end{equation}
and then via a skillful application of the maximum principle proved that these inequalities are preserved.

Finally, if $\{ M_t \}_{t\in (-\infty,T)}$ is an ancient solution, namely a solution that is defined for all $t\ll 0$, e.g. a blowup limit, then the maximum principle had infinite time to improve the evolving surface, which yields:

\begin{theorem}[noncollapsing and convexity \cite{HK_inscribed,HK1}]\label{thm_HK}
Every ancient $\alpha$-noncollapsed flow is in fact $1$-noncollapsed and convex.
\end{theorem}

Closely related convexity estimates have been obtained earlier by White \cite{white_nature} and Huisken-Sinestrari \cite{huisken-sinestrari1,huisken-sinestrari2}, and a sharp noncollapsing estimate has been obtained first by Brendle \cite{Brendle_inscribed}.

\subsection{Monotonicity formula and epsilon-regularity}

Recall that by equation \eqref{eq_areamon} the total area is monotone under mean curvature flow. However, since $\textrm{Area}(\lambda M)=\lambda^2 \textrm{Area}(M)$, this is not that useful when considering blowup sequences with $\lambda\to \infty$. A great advance was made by Huisken, who discovered a scale invariant monotone quantity.
To describe this, let $X_0=(x_0,t_0)$ be a point in space-time, and let
\begin{equation}
 \rho_{X_0}(x,t)=\frac{1}{4\pi(t_0-t)} e^{-\frac{|x-x_0|^2}{4(t_0-t)}}
\end{equation}
be the $2$-dimensional adjoint heat kernel centered at $X_0$, where $t<t_0$.

\begin{theorem}[Huisken's monotonicity formula \cite{Huisken_monotonicity}]\label{thm_huisken_mon} Along mean curvature flow we have
\begin{equation}\label{eq_huisken_mon}
 \frac{d}{dt}\int_{M_t} \rho_{X_0} d\mu = -\int_{M_t} \left|\vec{H}-\frac{(x-x_0)^\perp}{2(t-t_0)}\right|^2 \rho_{X_0} d\mu.
\end{equation}
\end{theorem} 

Huisken's monotonicity formula \eqref{eq_huisken_mon} can be thought of as weighted version of \eqref{eq_areamon}. A key property is its invariance under the parabolic rescaling $(x,t)\mapsto (\lambda (x-x_0),\lambda^2(t-t_0))$. Another key property is that the equality case of \eqref{eq_huisken_mon} exactly characterizes the selfsimilarly shrinking solutions. Indeed, centering at $(x_0,t_0)=(0,0)$ for ease of notation, note that for an ancient mean curvature flow $\{M_t\subset \mathbb{R}^{3}\}_{t\in (-\infty,0)}$ we have
\begin{equation}
\vec{H}-\frac{x^\perp}{2t}=0\quad \forall t<0  \qquad \Leftrightarrow \qquad M_t=\sqrt{-t}M_{-1}\quad \forall t<0.
\end{equation}
The monotone quantity appearing on the left hand side of \eqref{eq_huisken_mon},
\begin{equation}
\Theta(\mathcal{M},X_0,r):=\int_{M_{t_0-r^2}} \rho_{X_0} d\mu,
\end{equation}
is called the Gaussian density. Note that $\Theta(\mathcal{M},X_0,r)\equiv 1$ for all $r>0$ if and only if $\mathcal{M}$ is a multiplicity-one plane containing $X_0$. This can be upgraded to the following $\varepsilon$-regularity theorem:

\begin{theorem}[$\varepsilon$-regularity \cite{brakke,white_regularity,KasaiTonegawa,DPGS}]\label{app_thm_easy_brakke}
There exist universal constants $\varepsilon>0$ and $C<\infty$, such that
\begin{equation}
 \sup_{X\in P(X_0,r)}\Theta(\mathcal{M},X,r)<1+\varepsilon\quad \Rightarrow\quad
 \sup_{P(X_0,r/2)} |A|\leq {C}r^{-1}.
\end{equation}
\end{theorem}

Here, because of the parabolic scaling, the estimate naturally takes place in the parabolic ball
\begin{equation}
P(X_0,r)=B(x_0,r)\times (t_0-r^2,t_0].
\end{equation}

\subsection{Weak solutions}

We will now discuss weak solutions, specifically level set and Brakke solutions, which allow one to continue the evolution through singularities.

Level set flows have been studied originally in the framework of viscosity solutions for degenerate PDEs  \cite{OsherSethian,evans-spruck,CGG}. Later, Ilmanen found a completely elementary geometry reformulation \cite{Ilmanen_set}. To describe this, recall that by the avoidance principle smooth mean curvature flows do not bump into each other. Motivated by this, a family of closed sets $\{C_t\}_{t\geq 0}$ is called a \emph{subsolution} if it avoids all smooth solutions, namely
\begin{equation}\label{eq_avoid}
C_{t_0}\cap M_{t_0}=\emptyset \qquad \Rightarrow \qquad C_{t}\cap M_{t}=\emptyset \quad \forall t\in[t_0,t_1]
\end{equation}
for any smooth closed mean curvature flow $\{M_t\}_{t\in[t_0,t_1]}$.

\begin{definition}[level set flow] The \emph{level set flow} $\{F_t(C)\}_{t\geq 0}$ of any closed set $C$ is the maximal subsolution $\{C_t\}_{t\geq 0}$ with $C_0=C$.
\end{definition}

Considering the closure of the union of all subsolutions, namely
\begin{equation}
F_{t}(C)=\overline{\bigcup \{ C_{t} \, | \ \{C_{t'}\}_{t'\geq 0} \textrm{ is a subsolution} \}}\, ,
\end{equation}
one sees that the level set flow exists and is unique. Moreover, uniqueness readily implies the following basic properties:
\begin{itemize}
\item semigroup property: $F_0(C)=C$, $F_{t+t'}(C)=F_t(F_{t'}(C))$.
\item commutes with translations:  $F_t(C+x)=F_t(C)+x$.
\item containment: If $C\subseteq C'$, then $F_t(C)\subseteq F_t(C')$.
\end{itemize}

It is also not hard to see that level set solutions are consistent with classical solutions, namely if $\{M\}_{t\in [0,T)}$ is a smooth mean curvature flow of closed surfaces, then
\begin{equation}
F_t(M)=M_t\qquad \forall t\in[0,T).
\end{equation}
Furthermore, by interposing a $C^{1}$-surface one can check that level set flows also avoid each other \cite{Ilmanen_set,HershkovitsWhite_avoidance}, namely
\begin{equation}\label{eq_avoid_lev}
C\cap C' =\emptyset \qquad \Rightarrow\qquad F_t(C)\cap F_t(C')=\emptyset \quad \forall t\geq 0 ,
\end{equation}
provided that at least one of $C,C'$ is compact. While the level set solution is unique by definition, the evolution can be nonunique:

\begin{example}[fattening, \cite{White_ICM,IlmanenWhite,Ketover,LeeZhao,CDHS}]\label{ex_fattening}
There exists a closed embedded surface $M\subset\mathbb{R}^3$, that encounters a conical singularity after which $F_t(M)$ develops nonempty interior.
\end{example}

For the sake of intuition, it helps to compare this with the nonuniqueness/fattening of the figure "X" under curve shortening flow (see \cite{White_ICM} for a more detailed discussion).

To capture nonuniqueness more precisely, given any closed embedded surface $M\subset\mathbb{R}^3$, we denote by $K$ the compact domain enclosed by $M$, and define $\mathcal{K}$ as the space-time track of the level set flow of $K$, namely
\begin{equation}
\mathcal{K}:=\left\{ (x,t)\in\mathbb{R}^3\times [0,\infty)\, | \, x\in F_t(K)\right\} .
\end{equation}
Similarly, setting $K':=\overline{K^c}$, we observe that $\partial K' = \partial K=M$, and define $\mathcal{K}'$ as the space-time track of the level set flow of $K'$.

\begin{definition}[outer and inner flow \cite{Ilmanen_lectures,HershkovitsWhite_uniqueness}]\label{def_innerouter}
The \emph{outer flow} and \emph{inner flow} of $M$ are defined by
\begin{equation}
M_t:= \{ x \in \mathbb{R}^3 \, | \, (x,t)\in \partial\mathcal{K} \} \qquad\mathrm{and}\qquad
M_t':= \{ x \in \mathbb{R}^3 \, | \, (x,t)\in \partial\mathcal{K}' \} .
\end{equation}
\end{definition}
Here, it is most convenient to work with the boundary of space-time sets, but alternatively one can check that
\begin{equation}
M_{t}=\lim_{t'\nearrow t} \partial F_{t'}(K) ,
\end{equation}
and similarly for the inner flow.

The following captures the first time when nonuniqueness happens:
\begin{definition}[discrepancy time \cite{HershkovitsWhite_uniqueness}] The \emph{discrepancy time} is 
\begin{equation}
T_{\mathrm{disc}}:= \inf\{\,  t>0\, | \, M_t\neq M_t'\, \} \in (0,\infty]\, .
\end{equation}
\end{definition}
Thanks to the work of Bamler-Kleiner \cite{BK_mult1}, the discrepancy time $T_{\textrm{disc}}$ is in fact equal to the fattening time
\begin{equation}
T_{\mathrm{fat}}:=  \inf\{\,  t>0\, | \, \textrm{Int}(F_t(M)) \neq 0 \}.
\end{equation}
However, for applications it is usually sufficient to use the elementary inequality $T_{\mathrm{fat}}\geq T_{\mathrm{disc}}$.\\

While level set solutions are very well suited for discussing the question of uniqueness versus nonuniqueness, we also need another notion of solutions, so called Brakke flows, which is better suited for arguments based on the monotonicity formula and for passing to limits. To discuss this, recall that a Radon measure $\mu$ in $\mathbb{R}^{3}$ is integer two-recifiable, if at almost every point it possess a tangent plane of integer multiplicity. Namely, considering the rescaled measure $\mu_{x,\lambda}(A):=\lambda^{-2}\mu(\lambda A+x)$, for $\mu$-a.e. $x$ we have
\begin{equation}
\lim_{\lambda\to 0} \mu_{x,\lambda} = \theta \mathcal{H}^2\lfloor P\, ,
\end{equation}
for some positive integer $\theta$ and some plane $P$. We write $P=T_x\mu$. Also recall that the associated integral varifold is defined by
\begin{equation}
V_{\mu}(\psi)= \int \psi(x,T_x\mu) \, d\mu(x)\, .
\end{equation}

\begin{definition}[Brakke flow \cite{brakke,Ilmanen}]\label{def_brakke}
A two-dimensional \emph{integral Brakke flow} in $\mathbb{R}^{3}$ is a family of Radon measures $\mathcal M = \{\mu_t\}_{t\in I}$ in $\mathbb{R}^{3}$ that is integer two-rectifiable for almost every time and satisfies
\begin{equation}\label{Brakke_inequ}
\frac{d}{dt} \int \varphi \, d\mu_t \leq \int \left(D\varphi \cdot \vec{H} -\varphi {\vec{H}}^2 \right)\, d\mu_t
\end{equation}
for all test functions $\varphi\in C^1_c(\mathbb{R}^{3},\mathbb{R}_+)$. Here, $\tfrac{d}{dt}$ denotes the limsup of difference quotients, and ${\vec{H}}$ denotes the mean curvature vector of the associated varifold $V_{\mu_t}$.
\end{definition}

Here, by convention the right hand side of \eqref{Brakke_inequ} is interpreted as $-\infty$ whenever it does not make sense. Hence, it actually does make sense at almost every time.
Definition \ref{def_brakke} is motivated by fact that for any smooth flow $M_t$, the area measure $\mu_t=\mathcal{H}^2\lfloor M_t$ satisfies
 \eqref{Brakke_inequ} with equality, but only the inequality is preserved under weak limits. Using the definition, one can easily check that Huisken's monotone quantity still satisfies
\begin{equation}\label{app_loc_mon_brakke}
 \frac{d}{dt} \int\rho_{X_0}d\mu_t \leq -\int \left|\vec{H}-\frac{(x-x_0)^\perp}{2(t-t_0)}\right|^2 \rho_{X_0}d\mu_t\, .
\end{equation}
The following theorem enables one to pass to limits of Brakke flows:

\begin{theorem}[compactness \cite{brakke,Ilmanen}]\label{thm_comp}
Any sequence of integral Brakke flows $\{\mu^i_t\}$  with uniform area bounds on compact subsets has a subsequence that converges to an integral Brakke flow $\{\mu_t\}$.
\end{theorem}

Here, the convergence is in the sense of Radon measures at every time and in the sense of varifolds at almost every time. Also, the area bound hypothesis is always satisfied in practice thanks to \eqref{app_loc_mon_brakke}.

All integral Brakke flows that we will encounter are
\emph{unit-regular} \cite{white_regularity,SchulzeWhite},  i.e. near every space-time point of Gaussian density $1$ the flow is regular in a two-sided parabolic ball, and
\emph{cyclic} \cite{White_cyclic},  i.e. for a.e. $t$ the associated $\mathbb{Z}_2$ flat chain $[V_{\mu_t}]$ satisfies $\partial [V_{\mu_t}]=0$.\footnote{Intuitively, this simply means that we can color the inside and outside.} By the cited references, both properties are preserved under weak limits.\\

Finally, the notions of outer/inner flow and Brakke flow are compatible, as observed by Hershkovits-White \cite{HershkovitsWhite_uniqueness}. Specifically, given any closed initial surface $M_0$, using Ilmanen's elliptic regularization \cite{Ilmanen}, one can construct unit-regular, cyclic, integral Brakke flows $\{\mu_t\}$ and $\{\mu_t'\}$ with initial condition $\mathcal{H}^2\lfloor M_0$, such that their supports are given by the outer flow $\{M_t\}$ and the inner flow $\{M_t'\}$, respectively.

\section{Mean curvature flow through singularities in $\mathbb{R}^3$}\label{section3}

In this section, we discuss mean curvature flow through singularities in $\mathbb{R}^3$. In the special case of mean-convex surfaces the flow through singularities has been understood already 25 years ago thanks to the deep work of White \cite{White_size,white_nature} (see also \cite{White_ICM} for a survey, and \cite{HK1} for a streamlined approach based on Andrews' noncollapsing result). In the general case there is now also a powerful and highly developed theory, thanks to significant developments over the last 10 years, including:
\begin{itemize}
\item the uniqueness result for cylindrical tangent flows
\item the classification of genus zero shrinkers
\item the proof of Ilmanen's mean-convex neighborhood conjecture
\item the proof of Huisken's genericity conjecture
\item the proof of Ilmanen's multiplicity-one conjecture
\end{itemize}

We will now discuss these results and their consequences in nonhistorical order. Throughout, we denote by $\mathcal{M}=\{M_t\}_{t\geq 0}$ the outer (or inner) flow starting at a closed embedded surface $M_0\subset\mathbb{R}^3$. (All arguments generalize quite easily to other ambient 3-manifolds, but for concreteness and ease of notation we stick to $\mathbb{R}^3$.)\\

The capture the basic structure of singularities, given any space-time point $X_0=(x_0,t_0)$ one considers a tangent-flow at $X_0$, namely
\begin{equation}
\hat{\mathcal{M}}_{X_0}:=\lim_{i'\to \infty}\mathcal{D}_{\lambda_i}(\mathcal{M}-X_0) ,
\end{equation}
i.e. one shifts $X_0$ to the space-time origin, parabolically dilates by $\lambda_i\to \infty$, and passes to a subsequential limit. Thanks to Husiken's monotonicity inequality \eqref{app_loc_mon_brakke} and the compactness theorem for integral Brakke flows (Theorem \ref{thm_comp}), tangent-flows always exit and are always self-similarly shrinking. However, the tangent flows provided by  the compactness theorem come with some integer multiplicity, which could a priori be bigger than one. This potential scenario was ruled out in a recent breakthrough by Bamler-Kleiner:

\begin{theorem}[multiplicity-one \cite{BK_mult1}]\label{thm_mult_one}
If $M_0\subset \mathbb{R}^3$ is a closed embedded surface, then all tangent flows have multiplicity one.
\end{theorem}

The main worrisome scenario to rule out is two sheets connected by catenoidal necks that result in a multiplicity-two blowup limit. To do so, the authors introduce a novel sheet separation function $\mathfrak{s}$ (loosely speaking, if the sheets are graphs of functions $u_1< u_2$, then $\mathfrak{s}\sim u_2-u_1$), and prove that away from the necks one has the crucial differential inequality
\begin{equation}
(\partial_t-\Delta)\log \mathfrak{s}\geq 0,
\end{equation}
which shows that to some extent the sheets are moving away from each other.
They then deal with the neck regions via delicate integral estimates, starting from an old observation of Ilmanen \cite{Ilmanen_mon}, who combined  Huisken's monotonicity formula and the Gauss-Bonnet theorem to derive $L^2$-bounds for the second fundamental form.

Theorem \ref{thm_mult_one} (multiplicity-one), together with Theorem \ref{app_thm_easy_brakke} ($\varepsilon$-regularity) and dimension reduction (see e.g. \cite{White_stratification,CHN}), immediately yields a sharp estimate for the size of the singular set $\mathcal{S}\subset \mathbb{R}^3\times \mathbb{R}_+$.

\begin{corollary}[partial regularity \cite{BK_mult1}]
For the (outer or inner) flow of embedded surfaces the singular set $\mathcal{S}$ has  dimension at most 1.
\end{corollary}

Here, dimension refers to the Hausdorff (or Minkowski) dimension with respect to the parabolic metric on space-time, where time has the units of length squared, so e.g. $\mathrm{dim}(\mathbb{R}^3\times \mathbb{R}_+)=5$.

Thanks to Theorem \ref{thm_mult_one} (multiplicity-one) all singularities are modelled by multiplicity-one shrinkers. It is know that the ends of noncompact shrinkers are either cylindrical or asymptotically conical \cite{Wang_shrinker,BK_mult1}. Mean curvature flow through conical singularities can be nonunique, as discussed in Example \ref{ex_fattening}. Fortunately, asymptotically conical shrinkers can be excluded most of the time thanks to the following two fundamental results.

\begin{theorem}[genus zero shrinkers \cite{Brendle_sphere}]\label{thm_genus_zero}
The only nontrivial shrinkers of genus zero are the round sphere and the round cylinder.
\end{theorem}

This result, due to Brendle, is proved by an ingenious application of the stability inequality and the maximum principle. 

\begin{theorem}[generic singularities \cite{CCMS,CCS,CCMS2}]
For generic $M_0$, the only tangent flows at singular points are the round shrinking spheres and cylinders.
\end{theorem}

The result, due to Chodosh-Choi-Mantoulidis-Schulze,\footnote{See also the pioneering work by Colding-Minicozzi \cite{CM_generic}, and the alternative approach by Sun-Xue \cite{SunXue1,SunXue2}.} is proved most easily via a one-sided deformation argument, which shows that for any unstable singularity one can perturb such that the Gaussian density drops by a definite amount. (Historically, their original argument based on classifying one-sided ancient flows was much more involved, but after the solution of the multiplicity-one conjecture they gave a simple short proof.)

Observing also that spherical singularities are isolated, we can thus assume from now on that $\mathcal{M}$ has a neck-singularity at $X_0$, i.e.
\begin{equation}
\hat{\mathcal{M}}_{X_0}=\big\{ \mathbb{R}\times S^1(\sqrt{2|t|}) \big\}_{t<0}.
\end{equation}
This implicitly makes use of the following foundational result:

\begin{theorem}[uniqueness of tangent flows \cite{CM_uniqueness}]\label{thm_cm_uniq}
Cylindrical tangent flows are independent of the choice of rescaling factors  $\lambda_i\to \infty$.
\end{theorem}

This result, which in particular gives uniqueness of the axis, has been established via a Lojasiewicz-Simon argument  by Colding-Minicozzi, who overcame major difficulties caused by the noncompactness. We note that their argument also gives uniqueness of cylindrical tangent-flows at $-\infty$, which will be important below.

Finally, the key for mean curvature flow through neck-singularities is the proof of the mean-convex neighborhood conjecture from our joint work with Choi and Hershkovits:

\begin{theorem}[mean-convex neighborhoods \cite{CHH}]\label{thm_mean_convex}
If $\mathcal{M}=\{M_t\}_{t\geq 0}$ has a neck-singularity at $(x_0,t_0)$, then there exists $\varepsilon=\varepsilon(x_0,t_0)>0$ such that $M_t\cap B_\varepsilon(x_0)$ is mean-convex for $|t-t_0|<\varepsilon$.
\end{theorem}

Here, a major difficulty was to rule out the potential scenario of a degenerate neck-pinch with a non-convex cap. Namely, to fully capture the singularity, one really wants to understand all limit flows,
\begin{equation}
\mathcal{M}^\infty:=\lim_{\lambda_i\to \infty}\mathcal{D}_{\lambda_i}(\mathcal{M}-X_i) ,
\end{equation}
where now $X_i\to X_0$ depends on $i$ (e.g. for the degenerate neck-pinch one chooses $X_i$ along the tip). While tangent-flows are always self-similarly shrinking by Huisken's monotonicity formula, limit flows could be much more general. A priori, choosing $\lambda_i\to \infty$ suitably, we only know that $\mathcal{M}^\infty$ is an ancient asymptotically cylindrical flow, namely an ancient, unit-regular, cyclic, integral Brakke flow $\mathcal{M}=\{ \mu_t\}_{t\in (-\infty,T_E)}$ whose tangent flow at $-\infty$ is a round shrinking cylinder, i.e.
\begin{equation}\label{eq_tang_}
\lim_{\lambda\to 0} \mathcal{D}_{\lambda}\mathcal{M} = \{\mathbb{R}\times S^1(\sqrt{2|t|}) \}_{t<0}\, .
\end{equation}
Together with Choi and Hershkovits we classified all such flows:

\begin{theorem}[classification \cite{CHH}]\label{classification_theorem}
Any ancient asymptotically cylindrical flow in $\mathbb{R}^3$ is either a round shrinking cylinder, or a translating bowl, or an ancient oval.
\end{theorem}

\begin{figure}[htbp]
  \centering
\includegraphics[width=0.7\columnwidth]{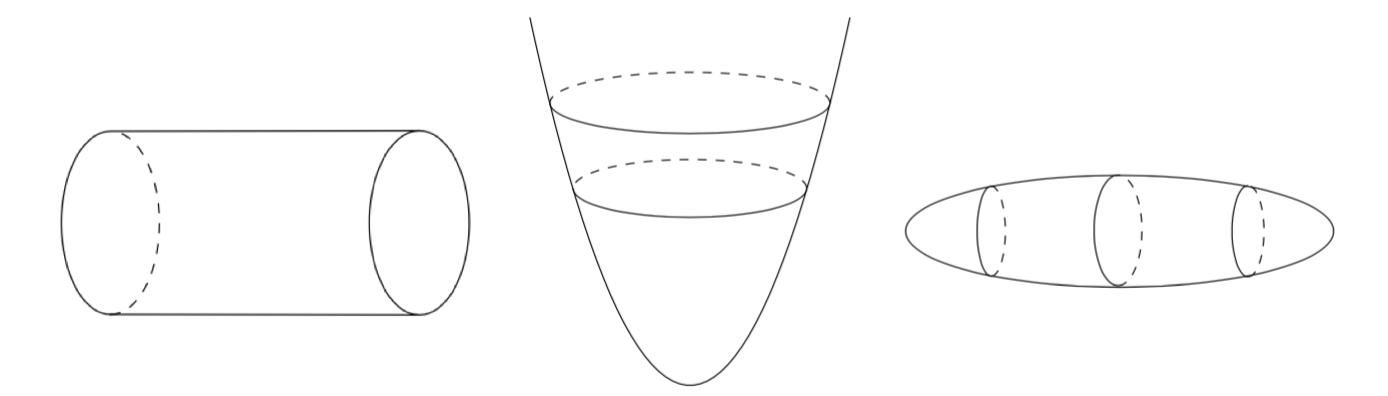}
\caption{Cylinder, bowl and ancient oval.}\label{class_asympt_cyl}
\end{figure}

The classification is illustrated in Figure \ref{class_asympt_cyl}. In particular, note that all solutions appearing in the classification result are convex. Hence, the mean-convex neighborhood conjecture follows from the classification theorem via a short argument by contradiction. Important prior classification results (under additional assumptions, such as convexity and/or selfsimilarity) can be found in \cite{Wang_convex,Haslhofer_bowl,BW,Hershkovits_translators}, and most importantly in the groundbreaking papers by Angenent-Daskalopoulos-Sesum \cite{ADS1,ADS2} and Brendle-Choi \cite{BC1}, which have been instrumental for our approach.

As an application of Theorem \ref{thm_mean_convex} (mean-convex neighborhoods), using also Theorem \ref{thm_mult_one} (multiplicity-one), one can confirm the uniqueness conjecture for mean curvature flow through neck-singularities:

\begin{corollary}[uniqueness \cite{HershkovitsWhite_uniqueness,CHH,BK_mult1}]
Mean curvature flow through neck-singularities is unique.
\end{corollary}

Indeed, it has been known since the 90s that nonfattening holds for mean-convex mean curvature flows \cite{evans-spruck,CGG}. More recently, Hershkovits-White \cite{HershkovitsWhite_uniqueness} localized this result and showed that it is enough to assume that all singularities have a mean-convex neighborhood, and we established exactly this assumption.

As another application, taking also into account Brendle's classification of genus zero shrinkers (Theorem \ref{thm_genus_zero}), one can confirm a conjecture of White, which he stated in his 2002 ICM lecture \cite{White_ICM}.

\begin{corollary}[flow of embedded 2-spheres \cite{Brendle_sphere,HershkovitsWhite_uniqueness,CHH,BK_mult1}]\label{cor_2spheres}
Mean curvature flow of embedded 2-spheres is well posed.
\end{corollary}

To conclude this section, let us outline the main steps of our proof of the classification theorem (Theorem \ref{classification_theorem}). Given any ancient asymptotically cylindrical flow $\mathcal{M}$ that is not a round shrinking cylinder we have to show that it is either a translating bowl or an ancient oval.

To get started, inspired by \cite{ADS1,BC1}, we set up a fine-neck analysis as follows. Given any $X_0=(x_0,t_0)$ we consider the renormalized flow
\begin{equation}
\bar{M}_\tau^{X_0}=e^{\tau/2}(M_{t_0-e^{-\tau}}-x_0) .
\end{equation}
Thanks to \eqref{eq_tang_} the flow $\bar{M}_\tau^{X_0}$ converges for $\tau\to -\infty$ to $\mathbb{R}\times S^1(\sqrt{2})$. Hence, writing $\bar{M}_\tau^{X_0}$ locally as a graph of a function $u^{X_0}(y,\vartheta,\tau)$ over the cylinder, the evolution is governed by the Ornstein-Uhlenbeck type operator
\begin{equation}
\mathcal{L}=\partial^2_y -\tfrac12 y \partial_y +\tfrac12 \partial_\vartheta^2+1 ,
\end{equation}
which is a self-adjoint operator on the Gaussian $L^2$-space. The operator $\mathcal{L}$ has 4 unstable eigenfunctions, namely $1$, $\sin \vartheta$, $\cos\vartheta$, $y$, and 3 neutral eigenfunctions, namely $y\sin\vartheta$, $y\cos\vartheta$, $y^2-2$, and all other eigenfunctions are stable. By the Merle-Zaag ODE lemma \cite{MZ} for $\tau \to -\infty$ either the unstable or neutral eigenfunctions dominate.

If the unstable case, we prove that there exists a constant $a=a(\mathcal{M})\neq 0$ independent of the center point $X_0$, such that (after suitable recentering to kill the rotations) we have
\begin{equation}\label{eq_fine_neck}
u^{X_0}(y,\vartheta,\tau)=aye^{\tau/2}+o(e^{\tau/2})
\end{equation}
for all $\tau\ll 0$ depending only on the cylindrical scale of $X_0$. Moreover, we show that every point outside a ball of controlled size in fact lies on such a fine-neck. Hence, by the Brendle-Choi neck-improvement theorem \cite{BC1} the solution becomes very symmetric at infinity. Finally, we can apply a variant of the moving plane method to conclude that the solution is smooth and rotationally symmetric, and hence a bowl. (For the last step, see also our joint paper with White \cite{CHHW}, where we developed a variant of the moving plane method for varifolds, in particular a Hopf lemma for Brakke flows, which allows us to derive symmetry and smoothness in tandem.)

Finally, in the neutral case, analyzing to ODE for the coefficient of the eigenfunction $y^2-2$, we show that there is an inwards quadratic bending, and consequently that the solution is compact. Blowing up near the tips, and taking into account that noncompact solutions have already been classified in the above paragraph, we see bowls. Hence, using the maximum principle we can show that the flow is mean-convex and noncollapsed. Therefore, we can apply the result by Angenent-Daskalopoulos-Sesum \cite{ADS2} to conclude that $\mathcal{M}$ is an ancient oval.

\bigskip

\section{Classification of noncollapsed singularities in $\mathbb{R}^4$}\label{section4}

In this section, we describe our recent classification of all noncollapsed singularities for the mean curvature flow in $\mathbb{R}^4$. We will assume throughout that all our ancient flows are noncollapsed in the sense of Definition \ref{def_andrews_static}, though we expect that this assumption eventually can be dropped in favour of the weaker assumption that the tangent flow at $-\infty$ is a round shrinking $\mathbb{R}^j\times S^{3-j}$. We also recall that ancient noncollapsed flows are always convex thanks to Theorem \ref{thm_HK}.

\begin{figure}[htbp]
  \centering
\includegraphics[width=0.5\columnwidth]{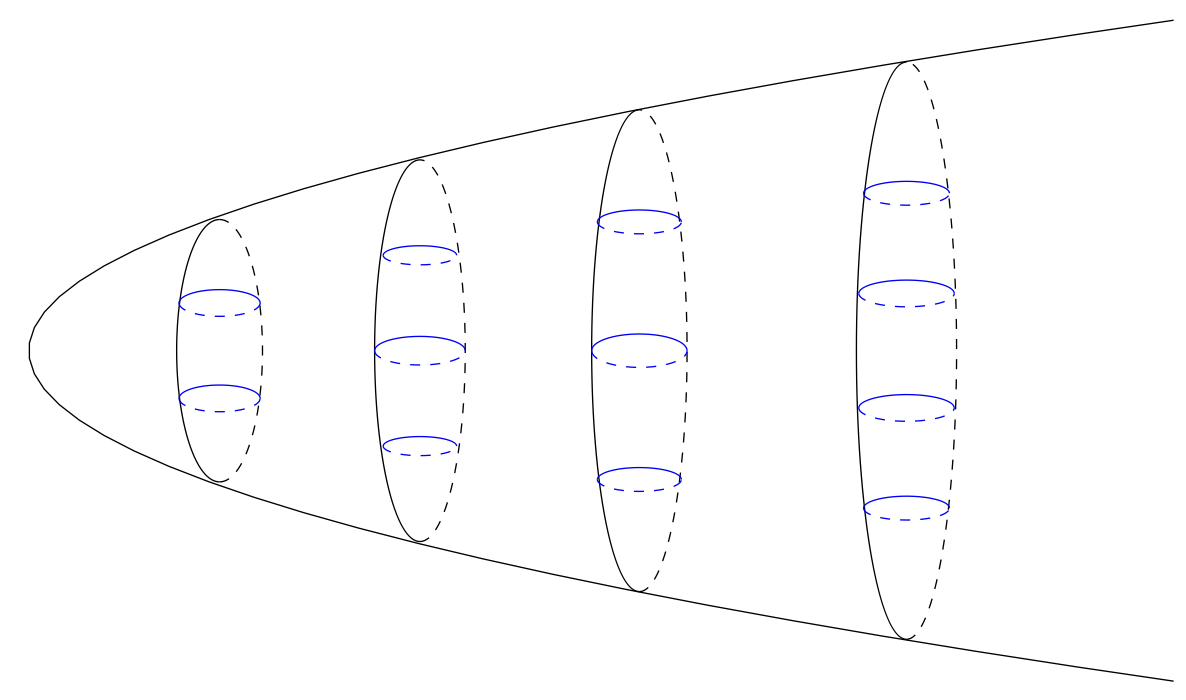}
\caption{The HIMW-translators look like oval-bowls.}\label{figure_oval_bowls}
\end{figure}

In stark contrast to $\mathbb{R}^3$, a classification of (noncollapsed) singularities in $\mathbb{R}^4$ until recently seemed out of reach. Fundamentally, this is because of the existence of examples with reduced symmetry. In fact, it is known since the pioneering work of Wang \cite{Wang_convex} that the Bernstein (i.e. symmetry) theorem fails for flows in $\mathbb{R}^4$. More precisely, Hoffman-Ilmanen-Martin-White constructed a 1-parameter family of 3d-translators that interpolate between the round 3d-bowl and $\mathbb{R}\times$2d-bowl, and are only $\mathbb{Z}_2\times \mathrm{O}_2$-symmetric \cite{HIMW}. These translators are illustrated in Figure \ref{figure_oval_bowls}.
Similarly, in joint work with Du we constructed a 1-parameter family of 3d-ovals that interpolate between the $\mathrm{O}_2\times \mathrm{O}_2$-symmetric 3d-oval and $\mathbb{R}\times$2d-oval, and are only $\mathbb{Z}_2^2\times \mathrm{O}_2$-symmetric \cite{DH_ovals}.

Recently, in a series of 6 papers joint with Choi, Choi, Daskalopoulos, Du, Hershkovits, and Sesum we obtained a complete classification:

\begin{theorem}[classification \cite{CHH_wing,CHH_translator,DH_shape,DH_no_rotation,CDDHS,CH_R4}]\label{thm_main}
Any ancient noncollapsed flow in $\mathbb{R}^4$ is, up to scaling and rigid motion,
\begin{itemize}
\item either one of the standard shrinkers $S^3$, $\mathbb{R}\times S^2$, $\mathbb{R}^2\times S^1$ or $\mathbb{R}^3$,
\item or the 3d-bowl, or $\mathbb{R}\times$2d-bowl, or belongs to the 1-parameter family of $\mathbb{Z}_2\times \mathrm{O}_2$-symmetric translators from \cite{HIMW},
\item or the $\mathbb{Z}_2\times \mathrm{O}_3$-symmetric 3d-oval, or the $\mathrm{O}_2\times \mathrm{O}_2$-symmetric 3d-oval, or $\mathbb{R}\times$2d-oval,
or belongs to the one-parameter family of $\mathbb{Z}_2^2\times \mathrm{O}_2$-symmetric 3d-ovals from \cite{DH_ovals}.
\end{itemize}
\end{theorem}

In addition to the 1-parameter families of translators and ovals discussed above, our list of course also contains all classical historical examples, in particular the 3d-bowl and the two examples of cohomogeneity-one 3d-ovals from \cite{white_nature} and \cite{HaslhoferHershkovits_ancient}, respectively. As an immediate consequence we obtain a complete classification of all potential blowup limits (and thus a canonical neighborhood theorem) for mean-convex flows in $\mathbb{R}^4$:

\begin{corollary}[blowup limits and canonical neighborhoods \cite{CHH_wing,CHH_translator,DH_shape,DH_no_rotation,CDDHS,CH_R4}]\label{cor_main}
For the mean curvature flow of mean-convex hypersurfaces in $\mathbb{R}^4$ (or in a 4-manifold) every blowup limit is given by one of the solutions from the above list. In particular, for every $\varepsilon>0$ there is an $H_\varepsilon=H_\varepsilon(M_0)<\infty$, such that around any space-time point $(p,t)$ with $H(p,t)\geq H_\varepsilon$ the flow is $\varepsilon$-close (after rescaling) to one of the solutions from the above list.
\end{corollary}

Indeed, it is known since the work of White \cite{white_nature} that all blowup limits of mean-convex flows are ancient noncollapsed flows (see also \cite{HK1} for a simpler proof based on Theorem \ref{thm_andrews}). More generally, in light of Ilmanen's multiplicity-one and mean-convex neighborhood conjecture, which are theorems in $\mathbb{R}^3$ but still open in $\mathbb{R}^4$, the conclusion of the corollary is also expected to hold for blowup limits near any generic singularity.\\

To outline the our approach, let $\mathcal{M}=\{M_t\}$ be an ancient noncollapsed mean curvature flow in $\mathbb{R}^4$ that is neither a static plane nor a round shrinking sphere. By general theory \cite{CM_uniqueness,HK1} the tangent flow at $-\infty$ is either a neck or a bubble-sheet. Since the neck-case has already been dealt with in \cite{ADS2,BC2}, we can  assume that
\begin{equation}\label{eq_tang_2}
\lim_{\lambda\to 0} \mathcal{D}_{\lambda}\mathcal{M} = \big\{\mathbb{R}^2\times S^1(\sqrt{2|t|}) \big\}_{t<0}\, .
\end{equation}
The analysis of such ancient solutions starts by considering the bubble-sheet function $u=u(y,\vartheta,\tau)$, which captures the deviation of the renormalized flow $\bar{M}_\tau=e^{\tau/2} M_{-e^{-\tau}}$ from the round bubble-sheet $\mathbb{R}^2\times S^1(\sqrt{2})$. The evolution of $u$ is governed by the Ornstein-Uhlenbeck type operator
\begin{equation}\label{prel_OU}
\mathcal L=\partial_{y_1}^2+\partial_{y_2}^2-\tfrac{y_1}{2} \partial_{y_1}-\tfrac{y_2}{2} \partial_{y_2}+\tfrac{1}{2} \partial_{\vartheta}^2+1.
\end{equation}
The operator $\mathcal{L}$ has the unstable eigenfunctions $1,y_1,y_2,\cos \vartheta, \sin \vartheta$, 
and the neutral eigenfunctions $y_1^2-2,y_2^2-2,y_1y_2,y_1\cos\vartheta,y_1\sin\vartheta,y_2\cos\vartheta,y_2\sin\vartheta$. 
Based on these spectral properties, and taking also into that the $\vartheta$-dependence is tiny thanks to Zhu's symmetry improvement result \cite{Zhu}, in joint work with Du we proved:

\begin{theorem}[normal form \cite{DH_no_rotation,DH_shape}]\label{thm_norm_form}
For $\tau\to -\infty$, in suitable coordinates, in Gaussian $L^2$-norm we have
\begin{equation}\label{bubble_sheet_quant1.1}
u = O(e^{\tau/2})\qquad\mathrm{or}\qquad
u= \frac{4-y_1^2-y_2^2}{\sqrt{8}|\tau|}+o(|\tau|^{-1})\qquad\mathrm{or}\qquad 
u= \frac{2-y_2^2}{\sqrt{8}|\tau|}+o(|\tau|^{-1}).
\end{equation}
\end{theorem}
Accordingly, the classification problem can be split up into 3 cases, which we call the case of fast convergence, slow convergence, and mixed convergence, respectively.

The case of fast convergence, which is easiest, has been settled in a joint paper with Choi and Hershkovits:

\begin{theorem}[no wings \cite{CHH_wing}]\label{thm_no_wings}
There are no wing-like ancient noncollapsed flows in $\mathbb{R}^4$. In particular, if the convergence is fast, then $\mathcal{M}$ is either a round shrinking $\mathbb{R}^2\times S^1$ or a translating $\mathbb{R}\times$2d-bowl.
\end{theorem}

To prove this, in a spirit similar to \eqref{eq_fine_neck}, we showed that
\begin{equation}\label{fine_expansion_fst}
u^X=(a_1 y_1 + a_2 y_2)e^{\tau/2}+ o(e^{\tau/2})
\end{equation}
for all $\tau$ negative enough depending only on the bubble-sheet scale. Analyzing this expansion along potential different edges, we concluded that $\mathcal{M}$ in fact splits off a line (hence is not wing-like) and is selfsimilar.

The case of slow convergence has been settled in joint work with B. Choi, Daskalopoulos, Du and Sesum:

\begin{theorem}[bubble-sheet ovals \cite{CDDHS}]\label{thm_bubble_sheet_ovals}
If the convergence is slow, then $\mathcal{M}$ is either the $\mathrm{O}_2\times \mathrm{O}_2$-symmetric 3d-oval, or belongs to the one-parameter family of $\mathbb{Z}_2^2\times \mathrm{O}_2$-symmetric 3d-ovals from \cite{DH_ovals}.
\end{theorem}

Regarding the proof, let us just mention that \eqref{bubble_sheet_quant1.1} in the case of slow convergence means inwards quadratic bending in all directions, which yields that $M_t$ is compact with axes of length approximately $\sqrt{2 |t|\log |t|}$. Hence, up to numerous technical challenges, the problem turned out to be amenable to the techniques from \cite{ADS1,ADS2}.

We also recall that the special case of selfsimilarly translating solutions (with any rate of convergence) has been dealt with earlier:
\begin{theorem}[translators \cite{CHH_translator}]\label{thm_translators_intro}
If $\mathcal{M}$ is selfsimilarly translating, then it is either $\mathbb{R}\times$2d-bowl, or belongs to the 1-parameter family of $\mathbb{Z}_2\times \mathrm{O}_2$-symmetric entire translators from \cite{HIMW}.
\end{theorem}

Finally, in joint work with Choi, we settled the most difficult case of mixed convergence without any selfsimilarity assumption:

\begin{theorem}[mixed convergence \cite{CH_R4}]\label{thm_translators_intro}
If the convergence is mixed, then $\mathcal{M}$ is either $\mathbb{R}\times$2d-oval or is selfsimilarly translating.
\end{theorem}

Loosely speaking, to capture the (dauntingly small) slope in $y_1$-direction, we consider the derivative $u_1^X=\partial u^X/\partial y_1$, which kills the leading order dependence on $y_2$, and prove that
\begin{equation}
u_1^X = a e^{\tau/2}+ o (e^{\tau/2}).
\end{equation}
This differential neck theorem, which goes vastly beyond \eqref{fine_expansion_fst} can then be used to conclude that $\mathcal{M}$ is noncompact (hence there are no exotic ovals) and either splits off a line or is selfsimilarly translating.

\section{Further directions and open problems}\label{section5}

In this final section, we discuss some of the most important open problems and conjectures for mean curvature flow through singularities. Many of them are of course classical, c.f. Ilmanen's problem list \cite{Ilmanen_problems} and White's 2002 ICM lecture \cite{White_ICM}, but we also propose some new ones.

\subsection{Open problems for mean curvature flow in $\mathbb{R}^3$}

This subsection is concerned with the most traditional setting of 2-dimensional surfaces evolving by mean curvature flow in $\mathbb{R}^3$.\\

To begin with, a very important open problem about singularity models is the following question of Ilmanen:

\begin{problem}[rigidity of the cylinder]\label{q_no_cyl}
Is the round cylinder the only complete embedded shrinker with a cylindrical end?
\end{problem}

To answer this in the affirmative one in particular has to rule out the scenario of shrinkers of mixed-type, where some ends are cylindrical and some ends are conical. The question is quite delicate, since the answer becomes negative if the completeness or embeddedness assumption is dropped \cite{Wang_cylinder}. Together with prior work of Wang \cite{Wang_shrinker} (see also \cite{BK_mult1}) a positive resolution  would imply the nice dichotomy that all singularities are either of conical-type or of neck-type. Related to this, one of the most fundamental open problems is the uniqueness of tangents:

\begin{conjecture}[uniqueness of tangent-flows]\label{conj_tan}
Tangent-flows are independent of the choice of sequence of rescaling factors $\lambda_i\to \infty$.
\end{conjecture}

By Theorem \ref{thm_cm_uniq}, and results of Schulze \cite{Schulze_tangent} and Chodosh-Schulze \cite{ChodoshSchulze}, this would follow from an affirmative answer to Problem \ref{q_no_cyl}. Alternatively, one may also try to establish uniqueness directly. Also closely related, motivated by a corresponding conjecture of Perelman for 3d Ricci flow \cite{Perelman}, we have:

\begin{conjecture}[bounded diameter]
The intrinsic diameter stays uniformly bounded as one approaches the first singular time.
\end{conjecture}

By a result of Du \cite{Du} (see also \cite{GH_diameter}) this would follow from an affirmative answer to Problem \ref{q_no_cyl}. Alternatively, one may also try to establish the diameter bound directly, possibly after solving Conjecture \ref{conj_tan}.\\

Next, let me discuss several conjectures about the fine structure of the singular set. To be specific one can assume that the flow through singularities is either the outer flow or the inner flow (see Definition \ref{def_innerouter}).

\begin{conjecture}[isolation]\label{conj_isol}
All singularities are isolated unless an entire tube shrinks to a closed curve.
\end{conjecture}

The conjecture is motivated by the principle that solutions of the level set flow, while only twice-differentiable, to some extent behave like analytic functions \cite{CM_arrival}. Weaker formulations of the conjecture are that there are only finitely many singular times or that singularities are isolated generically (see also the recent work of Sun-Xue \cite{SunXue2}, who proved that nondegenerate neck-singularities are isolated). Another well-known related open problem is:

\begin{problem}[selfsimilarity of blowup limits]\label{prob_selfsim}
Are blowup limits always selfsimilar?
\end{problem}

Recall that by Theorem \ref{classification_theorem} among ancient asymptotically cylindrical flows only the ancient ovals are not selfsimilar.
In joint work with B. Choi and Hershkovits we proved that the ancient ovals occur as blowup limit if and only if there is an accumulation of spherical singularities \cite{CHH_ovals}. Yet another related open problem is:

\begin{problem}[shrinking tubes]
Which closed curves can arise as singular set of a mean curvature flow of embedded surfaces?
\end{problem}

The only known example is the marriage-ring, which shrinks to a round circle. By an impressive general result of Colding-Minicozzi \cite{CM_singular}, up to potentially countably many points, the singular set is contained in finitely many Lipschitz curves. However, this leaves open the problem of constructing nontrivial examples and the problem of ruling out gaps, which is again tightly related to Conjecture \ref{conj_isol} and Problem \ref{prob_selfsim}.\\

Let us now discuss a problem that specifically assumes that the surface is topologically a two-sphere:

\begin{problem}[MCF proof of Smale conjecture]\label{prob_smale}
Give a mean curvature flow proof of the Smale conjecture.
\end{problem}

There is a geometric proof by Hatcher \cite{Hatcher} and a Ricci flow proof by Bamler-Kleiner \cite{BamlerKleiner_Smale}, but it would be nice to have a direct proof by mean curvature flow.
In one of its many formulations, Smale's conjecture states that the space of embedded two-spheres in $\mathbb{R}^3$ is contractible. The $\pi_0$-part follows from \cite{BHH}, thanks to Theorem \ref{cor_2spheres} and \cite{DanielsHolgate}, but dealing with higher homotopy groups would require more work (the strategy from \cite{BamlerKleiner_Smale} seems useful).\\

Finally, let me mention two somewhat more open ended problems:

\begin{problem}[selection principle]
Is there a selection principle for flowing out of conical singularities?
\end{problem}

It has been proposed by Dirr-Luckhaus-Novaga \cite{DLN} and Yip \cite{Yip} that considering mean curvature flow with space-time white noise the physically relevant solutions will be selected in the vanishing noise limit. See also the related recent work by Chodosh, Daniels-Holgate, and Schulze \cite{CDHS}, who proved  that the level set flow of a surface with a conical singularity fattens if and only if the level set flow of the cone fattens.

\begin{problem}[immersed surfaces]
Develop a theory of weak solutions for the flow of immersed surfaces.
\end{problem}

Many of the methods described in this article crucially rely on embeddedness. Recently, Lambert and M\"ader-Baumdicker successfully extended some results to the Alexandrov-immersed setting \cite{LMB}. However, dealing with general immersed surfaces seems to require some fundamentally new ideas.\\

\subsection{Open problems and further directions in higher dimensions}

In this subsection, we discuss some open problems and further directions in higher dimensions.\\

First, most of the problems and conjectures discussed in the previous subsection (except for Problem \ref{prob_smale}, which is specifically about 2-spheres) of course have natural generalizations to higher dimensions. For a higher dimensional version of Conjecture \ref{conj_tan} (uniqueness of tangent-flows) one may impose the additional assumption of smoothness, since the question whether or not nonsmooth tangent-flows are unique seems rather daunting.\\

Second, regarding the 5 significant developments listed in Section \ref{section3}, the uniqueness of cylindrical tangent-flows by Colding-Minicozzi actually holds in all dimensions, and  Brendle's genus zero result is of course purely 2-dimensional. The other 3 are still open in $\mathbb{R}^d$ for $d\geq 4$, so let us state them here, starting with Ilmanen's multiplicity one-conjecture (as usual, beyond the first singular time we consider the outer or inner flow):

\begin{conjecture}[multiplicity-one]\label{conj_mult_one}
For the mean curvature flow of embedded hypersurfaces all blowup limits have multiplicity one.
\end{conjecture}

Since the proof in $\mathbb{R}^3$ by Bamler-Kleiner \cite{BK_mult1} relies on integral estimates that are purely 2-dimensional, fundamentally new ideas seem to be needed to deal with higher dimensions. A good preliminary problem to get started would be to rule out the potential blowup limit of two static planes with orthogonal intersection. Next, let us state Huisken's genericity conjecture:

\begin{conjecture}[genericity]
Mean curvature flow with generic initial data only encounters spherical and cylindrical singularities.
\end{conjecture}

Recent work by Chodosh-Mantoulidis-Schulze \cite{CMS_low2} shows that in $\mathbb{R}^4$ and $\mathbb{R}^5$ this would follow from a generic strong version of Conjecture \ref{conj_mult_one}. On the other hand, higher dimensions (except possibly $\mathbb{R}^6$, which might be amenable to their methods) seem to require a more general approach going beyond one-sided perturbations. The main conjecture about cylindrical singularities is of course Ilmanen's mean-convex neighborhood conjecture:

\begin{conjecture}[mean-convex neighborhoods]\label{conj_mean_conv_nbd}
Cylindrical singularities have a mean-convex neighborhood.
\end{conjecture}

For neck-singularities, i.e. cylinders with only one $\mathbb{R}$-factor, this has been shown in \cite{CHHW}, but the general case is still open.
Work in progress with Choi suggests that Theorem \ref{thm_main} (classification) can be generalized to ancient flows asymptotic to a bubble-sheet, without assuming noncollapsing a priori. If successful, this will prove Conjecture \ref{conj_mean_conv_nbd} in $\mathbb{R}^4$, and thus yield a satisfying theory of mean curvature flow through singularities in $\mathbb{R}^4$ modulo Conjecture \ref{conj_mult_one}. Moreover, I am quite optimistic that our arguments will ultimately generalize to all dimensions, though this will of course first require a generalization of Theorem \ref{thm_main} to higher dimensions:

\begin{problem}[ancient noncollapsed flows]\label{prob_anc_noncoll}
Classify ancient noncollapsed flows in all dimensions.
\end{problem}

Recently, Du-Zhu \cite{DZ} and Choi-Du-Zhu \cite{CDZ} made substantial progress on this problem. Assuming this classification problem can be solved, the other main problem is then:

\begin{problem}[convexity]\label{prob_conv}
Show that every ancient asymptotically cylindrical flow is convex.
\end{problem}

A solution of Problem \ref{prob_anc_noncoll} and Problem \ref{prob_conv} would yield a proof of Conjecture \ref{conj_mean_conv_nbd} in all dimensions.\\

Finally, Theorem \ref{thm_main} also motivates several conjectures for 4d Ricci flow, c.f. \cite{Haslhofer_4dRicci}. To state the first one, recall that a \emph{$\kappa$-solution} is an ancient solution of the Ricci flow, complete with bounded curvature on compact time intervals, that has nonnegative curvature operator, positive scalar curvature, and is $\kappa$-noncollapsed at all scales.

\begin{conjecture}[$\kappa$-solutions in 4d Ricci flow]\label{conj_kappa}
Any $\kappa$-solution in 4d Ricci flow is, up to scaling and finite quotients, given by one of the following solutions.

\begin{itemize}
\item  shrinkers: $S^4$, $\mathbb{C}P^2$, $S^2\times S^2$, $\mathbb{R}\times S^3$ or $\mathbb{R}^2\times S^2$.
\item steadies: 4d Bryant soliton, the 3d Bryant soliton times a line, or belongs to the 1-parameter family of $\mathbb{Z}_2\times\mathrm{O}_3$-symmetric steady solitons constructed by Lai \cite{Lai}.
\item ovals: the $\mathbb{Z}_2\times\mathrm{O}_4$-symmetric 4d oval from Perelman \cite{Perelman}, the 3d oval times a line, the $\mathrm{O}_2\times\mathrm{O}_3$-symmetric 4d oval from Buttsworth \cite{Butt}, or belongs to the 1-parameter family of $\mathbb{Z}_2^2\times\mathrm{O}_3$-symmetric ovals from \cite{Haslhofer_4dRicci}.
\end{itemize}
\end{conjecture}

By results of Brendle-Naff \cite{BN} and Brendle-Daskalopoulos-Naff-Sesum \cite{BDNS} the conjecture holds true in case the tangent flow at $-\infty$ is a neck. So the remaining problem is to deal with the bubble-sheet case. In this direction, Ma-Mahmoudian-Sesum recently made some initial progress on the classification of steady solitons \cite{MMS}. More boldly, motivated by Theorem \ref{classification_theorem}, it is tempting to generalize Conjecture \ref{conj_kappa} as follows:

\begin{conjecture}[ancient 4d Ricci flows]\label{conj_cyl}
Any ancient asymptotically cylindrical 4d Ricci flow is, up to scaling, either $\mathbb{R}\times S^3$ or $\mathbb{R}^2\times S^2$ or one of the steadies or ovals listed in Conjecture \ref{conj_kappa}.
\end{conjecture}

In this direction, in a very interesting recent paper \cite{Law}, Law proved that any steady Ricci solitons asymptotic to $\mathbb{R}\times S^3$ is in fact the Bryant soliton (see also the earlier paper by Zhao-Zhu \cite{ZhaoZhu}, which assumes positive curvature outside a compact set). On the other hand, there are also interesting differences between Ricci flow and mean curvature flow. Specifically, I believe that quotient-necks lead to nonuniqueness:

\begin{conjecture}[nonuniqueness]\label{conj_nonunique}
Let $k\geq 3$. Then there exist a Ricci flow on $S^1\times S^3/\mathbb{Z}_k$ that forms a quotient-neck singularity with tangent flow $\mathbb{R}\times S^3/\mathbb{Z}_k$, and continuous nonuniquely modelled either on the Bryant soliton modulo $\mathbb{Z}_k$ or on Appleton's cohomogeneity-one soliton on the line bundle $O_{\mathbb{CP}^1}(-k)$ from \cite{Appleton_soliton}.
\end{conjecture}

Bamler-Kleiner proved that 3d Ricci flow through singularities is unique \cite{BK_uniqueness}, and Angenent-Knopf \cite{AK_conical} showed that conical singularities can cause nonuniqueness in dimension 5 and higher. So a solution of Conjecture \ref{conj_nonunique} would settle the uniqueness versus nonuniqueness problem in the remaining critical dimension 4.

\bibliography{haslhofer_references}

\newcommand{\etalchar}[1]{$^{#1}$}
\newcommand{\noopsort}[1]{} \newcommand{\singleletter}[1]{#1}
\begin{thebibliography}{CCMS24b}

\bibitem[AAG95]{AltschulerAngenentGiga}
S.~Altschuler, S.~Angenent, and Y.~Giga.
\newblock Mean curvature flow through singularities for surfaces of rotation.
\newblock {\em J. Geom. Anal.}, 5(3):293--358, 1995.

\bibitem[ADS19]{ADS1}
S.~Angenent, P.~Daskalopoulos, and N.~Sesum.
\newblock Unique asymptotics of ancient convex mean curvature flow solutions.
\newblock {\em J. Differential Geom.}, 111(3):381--455, 2019.

\bibitem[ADS20]{ADS2}
S.~Angenent, P.~Daskalopoulos, and N.~Sesum.
\newblock Uniqueness of two-convex closed ancient solutions to the mean
  curvature flow.
\newblock {\em Ann. of Math. (2)}, 192(2):353--436, 2020.

\bibitem[AK22]{AK_conical}
S.~Angenent and D.~Knopf.
\newblock Ricci solitons, conical singularities, and nonuniqueness.
\newblock {\em Geom. Funct. Anal.}, 32(3):411--489, 2022.

\bibitem[And12]{andrews1}
B.~Andrews.
\newblock Noncollapsing in mean-convex mean curvature flow.
\newblock {\em Geom. Topol.}, 16(3):1413--1418, 2012.

\bibitem[App17]{Appleton_soliton}
A.~Appleton.
\newblock A family of non-collapsed steady {R}icci solitons in even dimensions
  greater or equal to four.
\newblock 2017.

\bibitem[AV97]{AngenentVelazquez}
S.~Angenent and J.~Vel\'azquez.
\newblock Degenerate neckpinches in mean curvature flow.
\newblock {\em J. Reine Angew. Math.}, 482:15--66, 1997.

\bibitem[BC19]{BC1}
S.~Brendle and K.~Choi.
\newblock Uniqueness of convex ancient solutions to mean curvature flow in
  {$\Bbb R^3$}.
\newblock {\em Invent. Math.}, 217(1):35--76, 2019.

\bibitem[BC21]{BC2}
S.~Brendle and K.~Choi.
\newblock Uniqueness of convex ancient solutions to mean curvature flow in
  higher dimensions.
\newblock {\em Geom. Topol.}, 25(5):2195--2234, 2021.

\bibitem[BDNS23]{BDNS}
S.~Brendle, P.~Daskalopoulos, K.~Naff, and N.~Sesum.
\newblock Uniqueness of compact ancient solutions to the higher-dimensional
  {R}icci flow.
\newblock {\em J. Reine Angew. Math.}, 795:85--138, 2023.

\bibitem[BHH21]{BHH}
R.~Buzano, R.~Haslhofer, and O.~Hershkovits.
\newblock The moduli space of two-convex embedded spheres.
\newblock {\em J. Differential Geom.}, 118(2):189--221, 2021.

\bibitem[BK19]{BamlerKleiner_Smale}
R.~Bamler and B.~Kleiner.
\newblock Ricci flow and contractibility of spaces of metrics.
\newblock 2019.

\bibitem[BK22]{BK_uniqueness}
R.~Bamler and B.~Kleiner.
\newblock Uniqueness and stability of {R}icci flow through singularities.
\newblock {\em Acta Math.}, 228(1):1--215, 2022.

\bibitem[BK23]{BK_mult1}
R.~Bamler and B.~Kleiner.
\newblock On the multiplicity one conjecture for mean curvature flows of
  surfaces.
\newblock 2023.

\bibitem[BN23]{BN}
S.~Brendle and K.~Naff.
\newblock Rotational symmetry of ancient solutions to the {R}icci flow in
  higher dimensions.
\newblock {\em Geom. Topol.}, 27(1):153--226, 2023.

\bibitem[Bra78]{brakke}
K.~Brakke.
\newblock {\em The motion of a surface by its mean curvature}, volume~20 of
  {\em Mathematical Notes}.
\newblock Princeton University Press, Princeton, N.J., 1978.

\bibitem[Bre15]{Brendle_inscribed}
S.~Brendle.
\newblock A sharp bound for the inscribed radius under mean curvature flow.
\newblock {\em Invent. Math.}, 202(1):217--237, 2015.

\bibitem[Bre16]{Brendle_sphere}
S.~Brendle.
\newblock Embedded self-similar shrinkers of genus 0.
\newblock {\em Ann. of Math. (2)}, 183(2):715--728, 2016.

\bibitem[But22]{Butt}
T.~Buttsworth.
\newblock {$SO(2)\times SO(3)$}-invariant {R}icci solitons and ancient flows on
  {$\Bbb S^4$}.
\newblock {\em J. Lond. Math. Soc. (2)}, 106(2):1098--1130, 2022.

\bibitem[BW17]{BW}
J.~Bernstein and L.~Wang.
\newblock A topological property of asymptotically conical self-shrinkers of
  small entropy.
\newblock {\em Duke Math. J.}, 166(3):403--435, 2017.

\bibitem[CCMS24a]{CCMS}
O.~Chodosh, K.~Choi, C.~Mantoulidis, and F.~Schulze.
\newblock Mean curvature flow with generic initial data.
\newblock {\em Invent. Math.}, 237(1):121--220, 2024.

\bibitem[CCMS24b]{CCMS2}
O.~Chodosh, K.~Choi, C.~Mantoulidis, and F.~Schulze.
\newblock Revisiting generic mean curvature flow in $\mathbb{R}^3$.
\newblock 2024.

\bibitem[CCS23]{CCS}
O.~Chodosh, K.~Choi, and F.~Schulze.
\newblock Mean curvature flow with generic initial data {II}.
\newblock 2023.

\bibitem[CDD{\etalchar{+}}25]{CDDHS}
B.~Choi, P.~Daskalopoulos, W.~Du, R.~Haslhofer, and N.~Sesum.
\newblock Classification of bubble-sheet ovals in {$\Bbb{R}^4$}.
\newblock {\em Geom. Topol.}, 29(2):931--1016, 2025.

\bibitem[CDHS24]{CDHS}
O.~Chodosh, J.~Daniels-Holgate, and F.~Schulze.
\newblock Mean curvature flow from conical singularities.
\newblock {\em Invent. Math.}, 238(3):1041--1066, 2024.

\bibitem[CDZ25]{CDZ}
B.~Choi, W.~Du, and J.~Zhu.
\newblock Rigidity of ancient ovals in higher dimensional mean curvature flow.
\newblock 2025.

\bibitem[CGG91]{CGG}
Y.G. Chen, Y.~Giga, and S.~Goto.
\newblock Uniqueness and existence of viscosity solutions of generalized mean
  curvature flow equations.
\newblock {\em J. Differential Geom.}, 33(3):749--786, 1991.

\bibitem[CH24]{CH_R4}
K.~Choi and R.~Haslhofer.
\newblock Classification of ancient noncollapsed flows in $\mathbb{R}^4$.
\newblock 2024.

\bibitem[CHH21]{CHH_ovals}
B.~Choi, R.~Haslhofer, and O.~Hershkovits.
\newblock A note on the selfsimilarity of limit flows.
\newblock {\em Proc. Amer. Math. Soc.}, 149(3):1239--1245, 2021.

\bibitem[CHH22]{CHH}
K.~Choi, R.~Haslhofer, and O.~Hershkovits.
\newblock Ancient low-entropy flows, mean-convex neighborhoods, and uniqueness.
\newblock {\em Acta Math.}, 228(2):217--301, 2022.

\bibitem[CHH23]{CHH_translator}
K.~Choi, R.~Haslhofer, and O.~Hershkovits.
\newblock Classification of noncollapsed translators in {$\Bbb R^4$}.
\newblock {\em Camb. J. Math.}, 11(3):563--698, 2023.

\bibitem[CHH24]{CHH_wing}
K.~Choi, R.~Haslhofer, and O.~Hershkovits.
\newblock A nonexistence result for wing-like mean curvature flows in
  {$\Bbb{R}^4$}.
\newblock {\em Geom. Topol.}, 28(7):3095--3134, 2024.

\bibitem[CHHW22]{CHHW}
K.~Choi, R.~Haslhofer, O.~Hershkovits, and B.~White.
\newblock Ancient asymptotically cylindrical flows and applications.
\newblock {\em Invent. Math.}, 229(1):139--241, 2022.

\bibitem[CHN13]{CHN}
J.~Cheeger, R.~Haslhofer, and A.~Naber.
\newblock Quantitative stratification and the regularity of mean curvature
  flow.
\newblock {\em Geom. Funct. Anal.}, 23(3):828--847, 2013.

\bibitem[CM12]{CM_generic}
T.~Colding and W.~Minicozzi.
\newblock Generic mean curvature flow {I}; generic singularities.
\newblock {\em Ann. of Math. (2)}, 175(2):755--833, 2012.

\bibitem[CM15]{CM_uniqueness}
T.~Colding and W.~Minicozzi.
\newblock Uniqueness of blowups and {L}ojasiewicz inequalities.
\newblock {\em Ann. of Math. (2)}, 182(1):221--285, 2015.

\bibitem[CM16a]{CM_arrival}
T.~Colding and W.~Minicozzi.
\newblock Differentiability of the arrival time.
\newblock {\em Comm. Pure Appl. Math.}, 69(12):2349--2363, 2016.

\bibitem[CM16b]{CM_singular}
T.~Colding and W.~Minicozzi.
\newblock The singular set of mean curvature flow with generic singularities.
\newblock {\em Invent. Math.}, 204(2):443--471, 2016.

\bibitem[CMS]{CMS_low2}
O.~Chodosh, C.~Mantoulidis, and F.~Schulze.
\newblock Mean curvature flow with generic low-entropy data {II}.
\newblock {\em Duke Math. J. (to appear)}.

\bibitem[CS21]{ChodoshSchulze}
O.~Chodosh and F.~Schulze.
\newblock Uniqueness of asymptotically conical tangent flows.
\newblock {\em Duke Math. J.}, 170(16):3601--3657, 2021.

\bibitem[DH]{DH_ovals}
W.~Du and R.~Haslhofer.
\newblock On uniqueness and nonuniqueness of ancient ovals.
\newblock {\em Amer. J. Math. (to appear)}.

\bibitem[DH22]{DanielsHolgate}
J.~Daniels-Holgate.
\newblock Approximation of mean curvature flow with generic singularities by
  smooth flows with surgery.
\newblock {\em Adv. Math.}, 410:Paper No. 108715, 42, 2022.

\bibitem[DH23]{DH_no_rotation}
W.~Du and R.~Haslhofer.
\newblock A nonexistence result for rotating mean curvature flows in {$\Bbb
  R^4$}.
\newblock {\em J. Reine Angew. Math.}, 802:275--285, 2023.

\bibitem[DH24]{DH_shape}
W.~Du and R.~Haslhofer.
\newblock Hearing the shape of ancient noncollapsed flows in {$\Bbb R^4$}.
\newblock {\em Comm. Pure Appl. Math.}, 77(1):543--582, 2024.

\bibitem[DLN01]{DLN}
N.~Dirr, S.~Luckhaus, and M.~Novaga.
\newblock A stochastic selection principle in case of fattening for curvature
  flow.
\newblock {\em Calc. Var. Partial Differential Equations}, 13(4):405--425,
  2001.

\bibitem[DPGS24]{DPGS}
Guido De~Philippis, Carlo Gasparetto, and Felix Schulze.
\newblock A short proof of {A}llard's and {B}rakke's regularity theorems.
\newblock {\em Int. Math. Res. Not. IMRN}, (9):7594--7613, 2024.

\bibitem[Du21]{Du}
W.~Du.
\newblock Bounded diameter under mean curvature flow.
\newblock {\em J. Geom. Anal.}, 31(11):11114--11138, 2021.

\bibitem[DZ]{DZ}
W.~Du and J.~Zhu.
\newblock Spectral quantization for ancient asymptotically cylindrical flows.
\newblock {\em Adv. Math. (to appear)}.

\bibitem[ES91]{evans-spruck}
L.~Evans and J.~Spruck.
\newblock Motion of level sets by mean curvature. {I}.
\newblock {\em J. Differential Geom.}, 33(3):635--681, 1991.

\bibitem[GH20]{GH_diameter}
P.~Gianniotis and R.~Haslhofer.
\newblock Diameter and curvature control under mean curvature flow.
\newblock {\em Amer. J. Math.}, 142(6):1877--1896, 2020.

\bibitem[Ham95]{Hamilton_survey}
R.~Hamilton.
\newblock The formation of singularities in the {R}icci flow.
\newblock In {\em Surveys in differential geometry, {V}ol. {II} ({C}ambridge,
  {MA}, 1993)}, pages 7--136. Int. Press, Cambridge, MA, 1995.

\bibitem[Has15]{Haslhofer_bowl}
R.~Haslhofer.
\newblock Uniqueness of the bowl soliton.
\newblock {\em Geom. Topol.}, 19(4):2393--2406, 2015.

\bibitem[Has24a]{Haslhofer_summer_school}
R.~Haslhofer.
\newblock Lectures on mean curvature flow of surfaces.
\newblock 2024.

\bibitem[Has24b]{Haslhofer_4dRicci}
R.~Haslhofer.
\newblock On $\kappa$-solutions and canonical neighborhoods in 4d {R}icci flow.
\newblock {\em J. Reine Angew. Math.}, 811:257--265, 2024.

\bibitem[Hat83]{Hatcher}
A.~Hatcher.
\newblock A proof of the {S}male conjecture, {${\rm Diff}(S^{3})\simeq {\rm
  O}(4)$}.
\newblock {\em Ann. of Math. (2)}, 117(3):553--607, 1983.

\bibitem[Her20]{Hershkovits_translators}
O.~Hershkovits.
\newblock Translators asymptotic to cylinders.
\newblock {\em J. Reine Angew. Math.}, 766:61--71, 2020.

\bibitem[HH16]{HaslhoferHershkovits_ancient}
R.~Haslhofer and O.~Hershkovits.
\newblock Ancient solutions of the mean curvature flow.
\newblock {\em Comm. Anal. Geom.}, 24(3):593--604, 2016.

\bibitem[HI01]{HuiskenIlmanen}
G.~Huisken and T.~Ilmanen.
\newblock The inverse mean curvature flow and the {R}iemannian {P}enrose
  inequality.
\newblock {\em J. Differential Geom.}, 59(3):353--437, 2001.

\bibitem[HIMW19]{HIMW}
D.~Hoffman, T.~Ilmanen, F.~Martin, and B.~White.
\newblock Graphical translators for mean curvature flow.
\newblock {\em Calc. Var. Partial Differential Equations}, 58(4):Paper No. 117,
  29, 2019.

\bibitem[HK]{HaslhoferKetover2}
R.~Haslhofer and D.~Ketover.
\newblock Free boundary minimal disks in convex balls.
\newblock {\em J. Reine Angew. Math. (to appear)}.

\bibitem[HK15]{HK_inscribed}
R.~Haslhofer and B.~Kleiner.
\newblock On {B}rendle's estimate for the inscribed radius under mean curvature
  flow.
\newblock {\em Int. Math. Res. Not.}, (15):6558--6561, 2015.

\bibitem[HK17]{HK1}
R.~Haslhofer and B.~Kleiner.
\newblock Mean curvature flow of mean convex hypersurfaces.
\newblock {\em Comm. Pure Appl. Math.}, 70(3):511--546, 2017.

\bibitem[HK19]{HaslhoferKetover}
R.~Haslhofer and D.~Ketover.
\newblock Minimal 2-spheres in 3-spheres.
\newblock {\em Duke Math. J.}, 168(10):1929--1975, 2019.

\bibitem[HP99]{HP}
G.~Huisken and A.~Polden.
\newblock Geometric evolution equations for hypersurfaces.
\newblock In {\em Calculus of variations and geometric evolution problems
  ({C}etraro, 1996)}, volume 1713 of {\em Lecture Notes in Math.}, pages
  45--84. Springer, Berlin, 1999.

\bibitem[HS99a]{huisken-sinestrari1}
G.~Huisken and C.~Sinestrari.
\newblock Mean curvature flow singularities for mean convex surfaces.
\newblock {\em Calc. Var. Partial Differential Equations}, 8(1):1--14,
  {\noopsort{a}}1999.

\bibitem[HS99b]{huisken-sinestrari2}
G.~Huisken and C.~Sinestrari.
\newblock Convexity estimates for mean curvature flow and singularities of mean
  convex surfaces.
\newblock {\em Acta Math.}, 183(1):45--70, {\noopsort{b}}1999.

\bibitem[HS09]{huisken-sinestrari3}
G.~Huisken and C.~Sinestrari.
\newblock Mean curvature flow with surgeries of two-convex hypersurfaces.
\newblock {\em Invent. Math.}, 175(1):137--221, 2009.

\bibitem[Hui84]{Huisken_convex}
G.~Huisken.
\newblock Flow by mean curvature of convex surfaces into spheres.
\newblock {\em J. Differential Geom.}, 20(1):237--266, 1984.

\bibitem[Hui90]{Huisken_monotonicity}
G.~Huisken.
\newblock Asymptotic behavior for singularities of the mean curvature flow.
\newblock {\em J. Differential Geom.}, 31(1):285--299, 1990.

\bibitem[HW20]{HershkovitsWhite_uniqueness}
O.~Hershkovits and B.~White.
\newblock Nonfattening of mean curvature flow at singularities of mean convex
  type.
\newblock {\em Comm. Pure Appl. Math.}, 73(3):558--580, 2020.

\bibitem[HW23]{HershkovitsWhite_avoidance}
O.~Hershkovits and B.~White.
\newblock Avoidance for set-theoretic solutions of mean-curvature-type flows.
\newblock {\em Comm. Anal. Geom.}, 31(1):31--67, 2023.

\bibitem[Ilm93]{Ilmanen_set}
T.~Ilmanen.
\newblock The level-set flow on a manifold.
\newblock In {\em Differential geometry: partial differential equations on
  manifolds ({L}os {A}ngeles, {CA}, 1990)}, volume~54 of {\em Proc. Sympos.
  Pure Math.}, pages 193--204. Amer. Math. Soc., Providence, RI, 1993.

\bibitem[Ilm94]{Ilmanen}
T.~Ilmanen.
\newblock Elliptic regularization and partial regularity for motion by mean
  curvature.
\newblock {\em Mem. Amer. Math. Soc.}, 108(520):x+90, 1994.

\bibitem[Ilm95]{Ilmanen_mon}
T.~Ilmanen.
\newblock Singularities of mean curvature flow of surfaces.
\newblock 1995.

\bibitem[Ilm98]{Ilmanen_lectures}
T.~Ilmanen.
\newblock Lectures on mean curvature flow and related equations.
\newblock 1998.

\bibitem[Ilm03]{Ilmanen_problems}
T.~Ilmanen.
\newblock Problems in mean curvature flow.
\newblock 2003.

\bibitem[IW25]{IlmanenWhite}
T.~Ilmanen and B.~White.
\newblock Fattening in mean curvature flow.
\newblock {\em Ars Inven. Anal.}, Paper No. 4, 32, 2025.

\bibitem[Ket24]{Ketover}
D.~Ketover.
\newblock Self-shrinkers whose asymptotic cones fatten.
\newblock 2024.

\bibitem[KT14]{KasaiTonegawa}
K.~Kasai and Y.~Tonegawa.
\newblock A general regularity theory for weak mean curvature flow.
\newblock {\em Calc. Var. Partial Differential Equations}, 50(1-2):1--68, 2014.

\bibitem[Lai24]{Lai}
Y.~Lai.
\newblock A family of 3d steady gradient solitons that are flying wings.
\newblock {\em J. Differential Geom.}, 126(1):297--328, 2024.

\bibitem[Law25]{Law}
M.~Law.
\newblock Uniqueness of asymptotically cylindrical steady gradient {R}icci
  solitons.
\newblock 2025.

\bibitem[LM23]{LiokumovichMaximo}
Y.~Liokumovich and D.~Maximo.
\newblock Waist inequality for 3-manifolds with positive scalar curvature.
\newblock {\em Perspectives in scalar curvature. {V}ol. 2}, pages 799--831,
  2023.

\bibitem[LMB24]{LMB}
B.~Lambert and E.~M\"ader-Baumdicker.
\newblock A note on {A}lexandrov immersed mean curvature flow.
\newblock {\em J. Geom. Anal.}, 34(9):Paper No. 268, 25, 2024.

\bibitem[LZ24]{LeeZhao}
T.~Lee and X.~Zhao.
\newblock Closed mean curvature flows with asymptotically conical
  singularities.
\newblock 2024.

\bibitem[MM12]{MM}
C.~Mantegazza and L.~Martinazzi.
\newblock A note on quasilinear parabolic equations on manifolds.
\newblock {\em Ann. Sc. Norm. Super. Pisa Cl. Sci. (5)}, 11(4):857--874, 2012.

\bibitem[MMS23]{MMS}
Z.~Ma, H.~Mahmoudian, and N.~Sesum.
\newblock Unique asymptotics of steady {R}icci solitons with symmetry.
\newblock 2023.

\bibitem[Mul56]{Mullins}
W.~Mullins.
\newblock Two-dimensional motion of idealized grain boundaries.
\newblock {\em J. Appl. Phys.}, 27:900--904, 1956.

\bibitem[MZ98]{MZ}
F.~Merle and H.~Zaag.
\newblock Optimal estimates for blowup rate and behavior for nonlinear heat
  equations.
\newblock {\em Comm. Pure Appl. Math.}, 51(2):139--196, 1998.

\bibitem[OS88]{OsherSethian}
S.~Osher and J.~Sethian.
\newblock Fronts propagating with curvature-dependent speed: algorithms based
  on {H}amilton-{J}acobi formulations.
\newblock {\em J. Comput. Phys.}, 79(1):12--49, 1988.

\bibitem[Per02]{Perelman}
G.~Perelman.
\newblock The entropy formula for the {R}icci flow and its geometric
  applications.
\newblock 2002.

\bibitem[Sch08]{Schulze_isoperimetric}
F.~Schulze.
\newblock Nonlinear evolution by mean curvature and isoperimetric inequalities.
\newblock {\em J. Differential Geom.}, 79(2):197--241, 2008.

\bibitem[Sch14]{Schulze_tangent}
F.~Schulze.
\newblock Uniqueness of compact tangent flows in mean curvature flow.
\newblock {\em J. Reine Angew. Math.}, 690:163--172, 2014.

\bibitem[SW09]{ShengWang}
W.~Sheng and X.~Wang.
\newblock Singularity profile in the mean curvature flow.
\newblock {\em Methods Appl. Anal.}, 16(2):139--155, 2009.

\bibitem[SW20]{SchulzeWhite}
F.~Schulze and B.~White.
\newblock A local regularity theorem for mean curvature flow with triple edges.
\newblock {\em J. Reine Angew. Math.}, 758:281--305, 2020.

\bibitem[SX21]{SunXue1}
A.~Sun and J.~Xue.
\newblock Initial perturbation of the mean curvature flow for closed limit
  shrinker.
\newblock 2021.

\bibitem[SX22]{SunXue2}
A.~Sun and J.~Xue.
\newblock Generic mean curvature flows with cylindrical singularities.
\newblock 2022.

\bibitem[Wan11]{Wang_convex}
X.~Wang.
\newblock Convex solutions to the mean curvature flow.
\newblock {\em Ann. of Math. (2)}, 173(3):1185--1239, 2011.

\bibitem[Wan16a]{Wang_shrinker}
L.~Wang.
\newblock Asymptotic structure of self-shrinkers.
\newblock 2016.

\bibitem[Wan16b]{Wang_cylinder}
L.~Wang.
\newblock Uniqueness of self-similar shrinkers with asymptotically cylindrical
  ends.
\newblock {\em J. Reine Angew. Math.}, 715:207--230, 2016.

\bibitem[Whi97]{White_stratification}
B.~White.
\newblock Stratification of minimal surfaces, mean curvature flows, and
  harmonic maps.
\newblock {\em J. Reine Angew. Math.}, 488:1--35, 1997.

\bibitem[Whi00]{White_size}
B.~White.
\newblock The size of the singular set in mean curvature flow of mean convex
  sets.
\newblock {\em J. Amer. Math. Soc.}, 13(3):665--695, 2000.

\bibitem[Whi02]{White_ICM}
B.~White.
\newblock Evolution of curves and surfaces by mean curvature.
\newblock In {\em Proceedings of the {I}nternational {C}ongress of
  {M}athematicians, {V}ol. {I} ({B}eijing, 2002)}, pages 525--538. Higher Ed.
  Press, Beijing, 2002.

\bibitem[Whi03]{white_nature}
B.~White.
\newblock The nature of singularities in mean curvature flow of mean-convex
  sets.
\newblock {\em J. Amer. Math. Soc.}, 16(1):123--138, 2003.

\bibitem[Whi05]{white_regularity}
B.~White.
\newblock A local regularity theorem for mean curvature flow.
\newblock {\em Ann. of Math. (2)}, 161(3):1487--1519, 2005.

\bibitem[Whi09]{White_cyclic}
B.~White.
\newblock Currents and flat chains associated to varifolds, with an application
  to mean curvature flow.
\newblock {\em Duke Math. J.}, 148(1):41--62, 2009.

\bibitem[WZ23]{WangZhou}
Z.~Wang and X.~Zhou.
\newblock Existence of four minimal spheres in {$S^3$} with a bumpy metric.
\newblock 2023.

\bibitem[Yip98]{Yip}
N.~Yip.
\newblock Stochastic motion by mean curvature.
\newblock {\em Arch. Rational Mech. Anal.}, 144(4):313--355, 1998.

\bibitem[Zhu22]{Zhu}
J.~Zhu.
\newblock {$SO(2)$} symmetry of the translating solitons of the mean curvature
  flow in {$\Bbb{R}^4$}.
\newblock {\em Ann. PDE}, 8(1):Paper No. 6, 40, 2022.

\bibitem[ZZ23]{ZhaoZhu}
Z.~Zhao and X.~Zhu.
\newblock 4d steady gradient {R}icci solitons with nonnegative curvature away
  from a compact set.
\newblock {\em arXiv:2310.12529}, 2023.

\end{thebibliography}

\bibliographystyle{alpha}

\bigskip

{\sc Department of Mathematics, University of Toronto, 40 St George Street, Toronto, ON  M5S 2E4, Canada}

\emph{E-mail:} roberth@math.toronto.edu

\end{document}